\documentclass{article}
\usepackage{amsmath,amssymb,verbatim}

\usepackage[dvips]{graphicx}
\usepackage{amsfonts}
\usepackage{amsthm}
\usepackage{epsfig}

\theoremstyle{plain}
\newtheorem{theorem}{Theorem}[section]
\newtheorem{lemma}[theorem]{Lemma}
\newtheorem{proposition}[theorem]{Proposition}
\newtheorem{corollary}[theorem]{Corollary}
\theoremstyle{definition}
\newtheorem{definition}[theorem]{Definition}

\newtheorem{remark}[theorem]{Remark}

\newtheorem{example}[theorem]{Example}

\newtheorem{notation}[theorem]{Notation}

\newtheorem{algorithm}[theorem]{Algorithm}

\newcommand{\xs}{x_1,\ldots,x_n}                
\newcommand{\vs}{v_1,\ldots,v_n}                
\newcommand{\Ms}{M_1,\ldots,M_q}                
\newcommand{\Fs}{F_1,\ldots,F_q}                

\newcommand{\height}{{\rm{height}}\ }           
\newcommand{\dimn}{ {\rm{dim}} \ }
\newcommand{\F}{{\mathcal{F}}}                  
\newcommand{\D}{\Delta}                         
\newcommand{\al}{\alpha}                         
\newcommand{\df}{\delta_{\mathcal{F}}}          
              
\newcommand{\st}{\ | \ }                        
\newcommand{\tuple}[1]{\langle #1 \rangle}      
\newcommand{\stuple}[1]{\{ #1 \}}               
\newcommand{\srmv}[1]{\setminus \{ #1\}}        
\newcommand{\rmf}[1]{\setminus \{#1\}}          
\newcommand{\void}[1]{}

\newcommand{\website}{\texttt{%
http://\bk www.dm.unipi.it/\bk$\sim$caboara/\bk Research/}}
\newcommand{\bk}{\discretionary{}{}{}}

\newcommand{\calG}{\mathcal{G}}
\newcommand{\calF}{\mathcal{F}}
\newcommand{\cocoa}{\mbox{\rm C\kern-.13em o\kern-.07 em C\kern-.13em o\kern-.15em A}} 
\newcommand{\cocoax}{\mbox{C\kern-.13em o\kern-.07 em C\kern-.13em o\kern-.15em A}} 
\newcommand{\cocoal}{\mbox{\rm C\kern-.13em o\kern-.07 em C\kern-.13em o\kern-.15emA\kern-.1em L}}
\newcommand{\OO}[1]{{O}\!\left(#1\right)}

\newcommand{\sm}{\setminus}

\newcommand{\residue}[3]{#1_{#2}^{#3}}

\newcommand{\seq}{\subseteq}    
\newcommand{\es}{\emptyset}             
\newcommand{\s}[1]{\{#1\}}              
\newcommand{\abs}[1]{|{#1}|}            
\newcommand{\p}{\tuple}                 
\newcommand{\such}{\,\,|\,\,}           

\renewcommand{\leq}{\leqslant}          
\renewcommand{\geq}{\geqslant}  

\date{}

\author{Massimo Caboara\thanks{Department of Mathematics, University
of Pisa, caboara@dm.unipi.it.} \and 
Sara Faridi\thanks{Department of
Mathematics and Statistics, Dalhousie University, Halifax, Canada, 
faridi@mathstat.dal.ca. 
Research supported by NSERC.} 
\and
 Peter
Selinger\thanks{Department of
Mathematics and Statistics, Dalhousie University, Halifax, Canada, 
selinger@mathstat.dal.ca. 
Research supported by NSERC.}}

\title{Simplicial cycles and the computation of simplicial trees}
\begin{document}

\maketitle
\begin{abstract}
  We generalize the concept of a cycle from graphs to simplicial
  complexes.  We show that a simplicial cycle is either a sequence of
  facets connected in the shape of a circle, or is a cone over such a
  structure. We show that a simplicial tree is a connected cycle-free
  simplicial complex, and use this characterization to produce an
  algorithm that checks in polynomial time whether a simplicial
  complex is a tree.  We also present an efficient algorithm for
  checking whether a simplicial complex is grafted, and therefore
  Cohen-Macaulay.
\end{abstract}

\section{Introduction} 

The main goal of this paper is to demonstrate that it is possible to
check, in polynomial time, if a monomial ideal is the facet ideal of a
simplicial tree.

Facet ideals were introduced in~\cite{F1} (generalizing results
in~\cite{Vi1} and~\cite{SVV} on edge ideals of graphs) as a method to
study square-free monomial ideals. The idea is to associate a
simplicial complex to a square-free monomial ideal, where each facet
(maximal face) of the complex is the collection of variables that
appear in a monomial in the minimal generating set of the ideal (see
Definition~\ref{d:facetideal}). The ideal will then be called the
``facet ideal'' of this simplicial complex. A special class of
simplicial complexes are called ``simplicial trees''
(Definition~\ref{d:tree}).  The definition of a simplicial tree is a
generalization of the concept of a graph-tree. Facet ideals of trees
have many properties; for example, they have normal and Cohen-Macaulay
Rees rings~\cite{F1}. Finding such classes of ideals is in general a
difficult problem. Simplicial trees also have strong
Cohen-Macaulay properties: their facet ideals are always sequentially
Cohen-Macaulay \cite{F2}, and one can determine under precisely what
combinatorial conditions on the simplicial tree the facet ideal is
Cohen-Macaulay \cite{F3}.  In~\cite{F4} it is shown that the theory is
not restricted to square-free monomial ideals; via polarization, one
can extend many properties of facet ideals to all monomial ideals. All
these properties, and many others, make simplicial trees 
useful from an algebraic point of view.

But how does one determine if a given square-free monomial ideal is
the facet ideal of a simplicial tree? In Section~\ref{s:tree-problem},
we give a characterization of trees that shows this can be done in
polynomial time.  This characterization is based on a careful study of
the structure of cycles in Section~\ref{s:cycles}. The study of
simplicial cycles is indeed interesting in its own right.  In graph
theory, the concepts of a tree and of a cycle are closely linked to
each other: a tree is a connected graph that does not contain a cycle,
and a cycle is a minimal graph that is not a tree. Generalizing to the
simplicial case, we use the latter property, together with the
existing definition of a simplicial tree, to define the concept of a
simplicial cycle. We then prove the remarkable fact that a 
simplicial cycle is either a sequence of facets
connected in the shape of a circle, or a cone over such a 
structure. This in turns yields an alternative characterization of
trees, given in Section~\ref{s:tree-problem}.

This result enables us to produce a polynomial time
algorithm to decide whether a given simplicial complex is a tree. The
algorithm itself is introduced in Section~\ref{s:tree-algorithm},
where the complexity and optimizations are also discussed.
Section~\ref{s:CM-properties} focuses on the algebraic properties of
facet ideals: in Section~\ref{s:grafting} we discuss a method of
adding generators to a square-free monomial ideal (or facets to the
corresponding complex) so that the resulting facet ideal is
Cohen-Macaulay. This method is called ``grafting'' a simplicial
complex. For simplicial trees, being grafted and being Cohen-Macaulay
are equivalent conditions~\cite{F3}. We then introduce an algorithm
that checks whether or not a given simplicial complex is grafted and
discuss its complexity.

\paragraph*{Implementations.} 

The algorithms described in this paper have first been coded using
\cocoal, the programming language of the \cocoa\ system (please see http:/\!/cocoa.dima.unige.it/).  
These prototypical implementations can be downloaded
from~\cite{implementation}.  Much more efficient (but less user
friendly) C{\tt ++} implementations have been developed for several versions of 
Algorithm~\ref{a:tree-algorithm} using the {\cocoa}Lib 
framework (http:/\!/cocoa.dima.unige.it/\bk cocoalib/).  The C{\tt ++} code is also available at the website 
\cite{implementation}.

\section{Simplicial complexes and trees}\label{s:first-section} 

We define the basic notions related to facet ideals. More details and
examples can be found in~\cite{F1,F3}.

\begin{definition}[Simplicial complex, facet]
A \emph{simplicial complex} $\D$ over a finite set of vertices $V$
is a collection of subsets of $V$, with the property that if $F \in
\D$ then all subsets of $F$ are also in $\D$. An element of $\D$ is
called a \emph{face} of $\D$, and the maximal faces are called
\emph{facets} of $\D$.
\end{definition}

Since we are usually only interested in the facets, rather than all
faces, of a simplicial complex, it will be convenient to work with the
following definition:

\begin{definition}[Facet complex]
  A \emph{facet complex} over a finite set of vertices $V$ is a set
  $\D$ of subsets of $V$, such that for all $F,G\in\D$, $F\seq G$
  implies $F=G$.  Each $F\in\D$ is called a \emph{facet} of $\D$.
\end{definition}

\begin{remark}[Equivalence of simplicial complexes and facet complexes]
  The set of facets of a simplicial complex forms a facet
  complex. Conversely, the set of subsets of the facets of a facet
  complex is a simplicial complex. This defines a one-to-one
  correspondence between simplicial complexes and facet complexes. 
  In this paper, we will work primarily with facet complexes.
\end{remark}

We define facet ideals, giving a one-to-one correspondence
between facet complexes (or, equivalently, simplicial complexes) and
square-free monomial ideals.

\begin{definition}[Facet ideal of a facet complex, facet complex of an ideal]\label{d:facetideal} \ 
\begin{itemize}
\item Let $\D$ be a facet complex over a vertex set $\s{\vs}$. Let $k$
be a field, and let $R=k[\xs]$ be the polynomial ring with
indeterminates $\xs$. The \emph{facet ideal of} $\D$ is defined to be
the ideal of $R$ generated by all the square-free monomials
$x_{i_1}\ldots x_{i_s}$, where $\{v_{i_1},\ldots, v_{i_s}\}$ is a
facet of $\D$. We denote the facet ideal of $\D$ by $\F(\D)$.

\item Let $I=(\Ms)$ be an ideal in the polynomial ring $k[\xs]$, where
$k$ is a field and $\Ms$ are square-free monomials in $\xs$ that form
a minimal set of generators for $I$.  The \emph{facet complex of} $I$
is defined to be $\df(I)=\s{\Fs}$, where for each $i$, $F_i=\{v_j \st
x_j|M_i, \ 1 \leq j \leq n \}$.
\end{itemize}

\end{definition}

 From now on, we often use $\xs$ to denote both the vertices of $\D$
 and the variables appearing in $\F(\D)$.  We also sometimes ease the
 notation by denoting facets by their corresponding monomials; for
 example, we write $xyz$ for the facet $\{x,y,z\}$.

 We now generalize some notions from graph theory to facet complexes. Note
 that a graph can be regarded as a special kind of facet complex,
 namely one in which each facet has cardinality 2.

\begin{definition}[Path, connected facet complex]\label{def:connected}
  Let $\D$ be 
  a facet complex.  A sequence of facets $F_1,\ldots,F_n$ is
  called a \emph{path} if for all $i=1,\ldots,n-1$, $F_i\cap
  F_{i+1}\neq\es$.  We say that two facets $F$ and $G$ are {\em
  connected} in $\D$ if there exists a path $F_1,\ldots,F_n$ with
  $F_1=F$ and $F_n=G$. Finally, we say that $\D$ is
  \emph{connected} if every pair of facets is connected.
\end{definition}

\begin{notation}
  If $F$, $G$ and $H$ are facets of $\D$, $H\leq_F G$ means that
$H\cap F\seq G\cap F$. The relation $\leq_F$ defines a preorder
(reflexive and transitive relation) on the facet set of $\D$.
\end{notation}

\begin{definition}[Leaf, joint]\label{d:leaf} Let $F$ be a
 facet of a facet complex $\D$.  Then $F$ is called a \emph{leaf} of
$\D$ if either $F$ is the only facet of $\D$, or else there exists
some $G \in \D \rmf{F}$ such that for all $H\in\D\rmf{F}$, we have $H
\leq_F G$. The facet $G$ above is called a \emph{joint} of the leaf
$F$ if $F \cap G \neq \emptyset$.
\end{definition}

It follows immediately from the definition that every leaf $F$
contains at least one {\em free vertex}, i.e., a vertex that belongs
to no other facet.

\begin{example}\label{example11} 
In the facet complex $\D=\s{xyz, yzu, uv}$,  $xyz$ and $uv$ are leaves, but
 $yzu$ is not a leaf.  Similarly, in $\D'=\s{xyu,xyz,xzv}$, the only
 leaves are $xyu$ and $xzv$.
\[
\D=\begin{tabular}{c}\epsfig{file=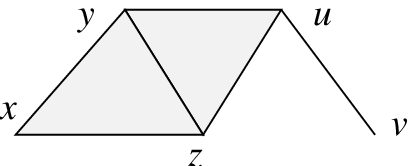, height=.6in}\end{tabular}
\qquad
\D'=\begin{tabular}{c}\epsfig{file=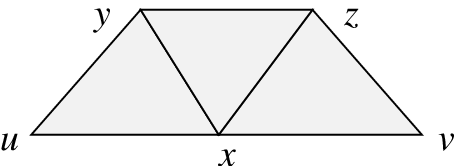, height=.6in}\end{tabular}
\]

\end{example}

\begin{definition}[Forest, tree]\label{d:tree} 
  A facet complex $\D$ is a \emph{forest} if every nonempty subset of
  $\D$ has a leaf. A connected forest is called a \emph{tree} (or
  sometimes a {\em simplicial tree} to distinguish it from a tree in
  the graph-theoretic sense).
\end{definition}

It is clear that any facet complex of cardinality one or two is a forest.
When $\D$ is a graph, the notion of a simplicial tree coincides with
that of a graph-theoretic tree.

\begin{example}\label{e:free-example} 
  The facet complexes in Example~\ref{example11} are trees. The facet
  complex pictured below has three leaves $F_1$, $F_2$ and $F_3$;
  however, it is not a tree, because if one removes the facet $F_4$,
  the remaining facet complex has no leaf.
\begin{center}
 \epsfig{file=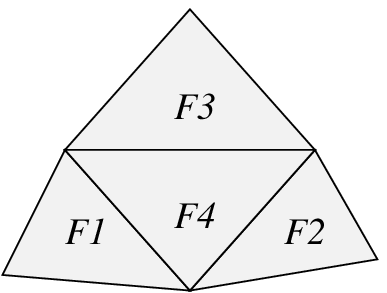, height=1in}
\end{center}
\end{example}

The following property is proved in {\cite[Lemma 4.1]{F3}}:

\begin{lemma}[A tree has two leaves]\label{l:twoleaves}
  Every tree with two or more facets has at least two leaves. \qed
\end{lemma}

\section{Cycles}\label{s:cycles}

In this section, we define a simplicial cycle as a minimal complex
without leaf. This in turns characterizes a simplicial tree as a
connected cycle-free facet complex. We further show that cycles
possess a particularly simple structure: each cycle is either
equivalent to a ``circle'' of facets with disjoint intersections, or
to a cone over such a circle. 

\begin{definition}[Cycle]\label{d:cycle}
A nonempty facet complex $\D$ is called a {\em cycle} if $\D$ has no
leaf but every nonempty proper subset of $\D$ has a leaf.
\end{definition}

Equivalently, $\D$ is a cycle if $\D$ is not a forest, but every
proper subset of $\D$ is a forest.  If $\D$ is a graph,
Definition~\ref{d:cycle} coincides with the graph-theoretic definition
of a cycle. The next two remarks are immediate consequences of the
definitions of cycle and forest:

\begin{remark}
  A cycle is connected.
\end{remark}

\begin{remark}
  A facet complex is a forest if and only if it does not contain a
  cycle.
\end{remark}

In the remainder of this section, we provide a complete
characterization of the structure of cycles.

\begin{definition}[Strong neighbor]\label{d:strong-neighbor} Let $\D$ 
be a facet complex and $F,G\in\D$. We say that $F$ and $G$ are
\emph{strong neighbors}, written $F\sim_{\D} G$, if $F\neq G$ and for
all $H\in\D$,
$F\cap G\seq H$ implies $H=F$ or $H=G$.
\end{definition}

  The relation $\sim_{\D}$ is symmetric, i.e., $F\sim_{\D} G$ if and only if
  $G\sim_{\D} F$. Note that if $\D$ has
  more than two facets, then $F\sim_{\D} G$ implies that $F \cap G \neq
  \emptyset$.

\begin{example} For the facet complex $\D'$ in Example~\ref{example11}, 
$xyu \not \sim_{\D'}xzv$, as their intersection $x$ lies in the facet
$xyz$. However, $xyz\sim_{\D'}xzv$ and similarly $xyz\sim_{\D'}xyu$.
\end{example}

\begin{remark}\label{r:strong-sub}
  Suppose $\D$ is a facet complex, and $\D'\seq\D$. Let
  $F,G\in\D'$. If $F\sim_{\D} G$, then $F\sim_{\D'} G$. The
  converse is not in general true.
\end{remark}

\begin{remark}\label{r:strong-leq} We have $F\sim_{\D} G$ if and only if  $G$ 
is strictly maximal with respect to $\leq_{F}$ on $\D\rmf{F}$, i.e.,
for all $H\neq F$, $G\leq_{F}H$ implies $G=H$. This is a simple
restatement of the definition.

\end{remark}

It turns out that a cycle can be described as a sequence of strong
neighbors. The following lemma follows directly from
Definition~\ref{d:strong-neighbor}.

\begin{lemma}\label{l:strong-double} If $\D$ is a facet complex with
 distinct facets $F,G_1,G_2$ such that $F\sim_{\D} G_1$ and $F\sim_{\D} G_2$,
then $F$ is not a leaf of $\D$.
\end{lemma}

    \begin{proof}
      If $F$ is a leaf, there exists a facet $H \neq F$
      such that $G_1\leq_F H$ and $G_2\leq_F H$, which by
      Remark~\ref{r:strong-leq} implies that $G_1=G_2=H$, a
      contradiction.
    \end{proof}

\begin{corollary}\label{c:cycle}
  Let $\D$ be a facet complex, and let $F_1,\ldots,F_n$ be distinct
facets with $n \geq 3$, such that $F_1\sim_{\D}F_2\sim_{\D} \ldots \sim_{\D} 
F_n\sim_{\D}F_1$. Then $\s{F_1,\ldots,F_n}$ has no leaf.
\end{corollary}
\begin{proof}
  This follows directly from Remark~\ref{r:strong-sub}, and
  Lemma~\ref{l:strong-double}.
\end{proof}

\begin{lemma}\label{l:subleaf}
  Suppose $\D$ is a facet complex and $F, G \in \D$. If $F$ is a leaf
  of $\D\sm\stuple{G}$, but not a leaf of $\D$, then $F\sim_{\D} G$.
\end{lemma}

\begin{proof}  Suppose $H$ is some facet such that
 $F\cap G\seq H$, but $H\neq F$ and $H\neq G$.  Since $F$ is a leaf
  for $\D\sm\stuple{G}$, there exists a facet $H'\in \D\sm\stuple{G}$
  such that $L\cap F\seq H'$ for all $L \in D\sm\stuple{F, G}$, and so
  $F\cap H\seq H'$.  But now we have $F\cap G\seq F\cap H\seq H'$,
  which implies that $F$ is a leaf of $\D$, a contradiction.
\end{proof}

\begin{proposition}[A cycle is a sequence of strong neighbors]\label{p:cycle}
  Suppose $\D$ is a cycle, and let $n=\abs{\D}$. Then $n\geq 3$, and
  the facets of $\D$ can be enumerated as $\D=\s{F_1,\ldots,F_n}$ in
  such a way that $$F_1\sim_{\D} F_2\sim_{\D} \ldots \sim_{\D}
  F_n\sim_{\D}F_1,$$ and $F_i\not\sim_{\D} F_j$ in all other cases, so that
  each facet is a strong neighbor of precisely two other facets.
\end{proposition}

    \begin{proof}
        First note that since $\D$ is not a forest, $n\geq 3$. We
        begin by showing that each facet has at least two distinct
        strong neighbors. Let $F\in\D$ be a facet. Since $\D$ is a
        cycle, $\D\sm\stuple{F}$ is a tree.  The subset
        $\D\sm\stuple{F}$ also has cardinality at least two, and
        therefore has two distinct leaves, say $G$ and $H$, by
        Lemma~\ref{l:twoleaves}. Since neither $G$ nor $H$ are leaves
        of $\D$ (because $\D$ is a cycle), we have $F\sim_{\D} G$ and
        $F\sim_{\D} H$ by Lemma~\ref{l:subleaf}.
  
        Now we can simply choose $F_1$ arbitrarily, then choose
        $F_2\neq F_1$ such that $F_1\sim_{\D} F_2$, then for every $i\geq
        3$ choose $F_i$ such that $F_{i-1}\sim_{\D} F_i$ and $F_i\neq
        F_{i-1},F_{i-2}$. Since $\D$ is finite, there will be some
        smallest $i$ such that $F_i=F_j$ for some $j<i$. Then
        $\D'=\stuple{F_j,\ldots,F_{i-1}}$ has no leaf by
        Corollary~\ref{c:cycle}, so $\D'=\D$. It follows that $j=1$
        and $i-1=n$. Finally, suppose that $F_k\sim_{\D} F_l$ for some
        $k\leq l-2$, where $k>1$ or $l<n$. Then
        $\stuple{F_1,\ldots,F_k,F_l,\ldots,F_n}$ has no leaf by
        Corollary~\ref{c:cycle}, contradicting the fact that it is a
        tree.
   \end{proof}

The converse of Proposition~\ref{p:cycle} is not true.

\begin{example}\label{e:cycle-is} The facet complex $\D$ is
 not a cycle, as its proper subset $\D'$ (which is indeed a cycle) has
no leaf. However, we have
$F_1\sim_{\D}F_2\sim_{\D}G\sim_{\D}F_3\sim_{\D}F_4\sim_{\D}F_1$, and
these are the only pairs of strong neighbors in $\D$.

\[
\D=\begin{tabular}{c}\epsfig{file=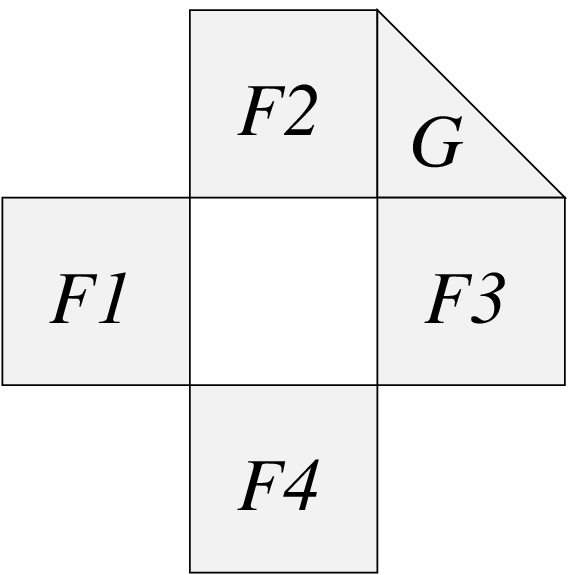, height=1in}\end{tabular}
\qquad
\D'=\begin{tabular}{c}\epsfig{file=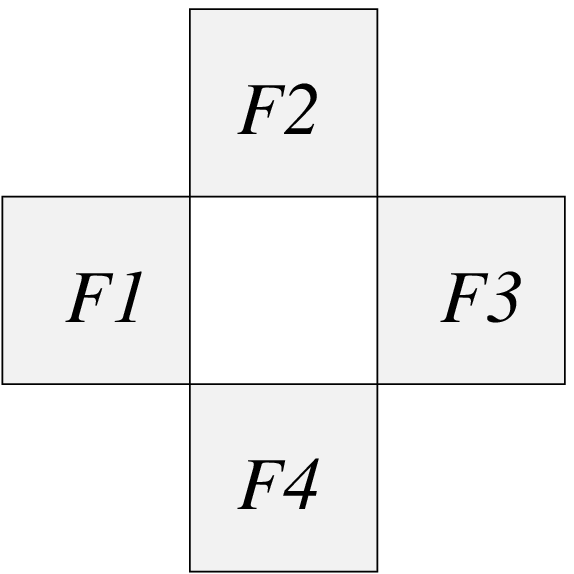, height=1in}\end{tabular}
\]

\end{example}

\begin{lemma}\label{l:cycle-to-tree} If $\D$ is a cycle, 
written as $F_1 \sim_{\D} F_2 \sim_{\D} \ldots \sim_{\D} F_n \sim_{\D}
 F_1$, then for each $i$, $\D_i=\D\srmv{F_i}$ is a tree with exactly
 two leaves $F_{i-1}$ and $F_{i+1}$, with joints $F_{i-2}$ and
 $F_{i+2}$, respectively.
\end{lemma}

    \begin{proof} We know that $\D_i$ is a tree, so it has 
      at least two leaves. By Lemma~\ref{l:strong-double} $F_{i-1}$
         and $F_{i+1}$ are the only choices. By
         Remark~\ref{r:strong-leq} $F_{i-2}$ is the only possible
         joint for $F_{i-1}$, and $F_{i+2}$ is the only possible joint
         for $F_{i+1}$.
    \end{proof}

The following lemma will be fundamental for the classification of cycles. 

\begin{lemma}\label{lem:cycle-fund}
  Let $\D$ be a cycle with facets $F\neq G\in\D$. If $F\not\sim_{\D}G$, then
  $F\cap G \seq H$ for all $H \in \D$.
\end{lemma}

\begin{proof}
  We first prove the claim in the special case where $F\sim_{\D}H$.
  Indeed, since $F$ is a strong neighbor of exactly two facets, there
  must be some $L\neq G,H$ such that $L\sim_{\D}F\sim_{\D}H$. Then
  Lemma~\ref{l:cycle-to-tree} implies that $H$ is a joint of $F$ in
  the tree $\D\srmv{L}$, and therefore $F\cap G \seq H$, or
  equivalently, $F\leq_G H$.
  
  Now consider the general case. By Proposition~\ref{p:cycle}, the
  facets of $\D$ can be enumerated as $F_1 \sim_{\D} F_2 \sim_{\D}
  \ldots \sim_{\D} F_n \sim_{\D} F_1$. Assume, without loss of
  generality, that $F=F_1$ and $G=F_i$, where $2<i<n$. By repeated
  applications of the special case above, we have
$$F\leq_G F_2  \leq_G \ldots \leq_G F_{i-1}.$$ 
In the other direction, we similarly have 
$$F\leq_G F_n \leq_G  F_{n-1} \leq_G  \ldots \leq_G F_{i+1}.$$
  Therefore, $F\cap G \seq F_j$ for $j=1,\ldots,n$.
\end{proof}

\begin{lemma}\label{l:cone} Let $\D$ be a facet complex, and let 
$$A=\bigcap_{F\in \D} F \mbox{ and } \D'=\{F\setminus A\st F \in \D\}.$$
Then $\D'$ is a facet complex. Moreover, $\D$ is a cycle if and only if $\D'$ is a cycle. \end{lemma}

         \begin{proof} 
           For each $F\in \D$, let $F'=F \setminus A$.  Since $\D$ is
           a facet complex, we have $F\not\seq G$ for any two distinct
           facets $F,G\in\D$, which clearly implies $F'\not\seq G'$.
           So $\D'$ is a facet complex. Let $\Gamma$ be any subset of
           $\D$, and let $\Gamma'=\s{F'\such F\in\D}$ be the
           corresponding subset of $\D'$. Then for any triple of
           facets $F,G,H \in \Gamma$, we have $F\leq_H G \iff F'
           \leq_{H'} G'$. Therefore, $\Gamma$ has a leaf if and only
           if $\Gamma'$ has a leaf. It follows that $\D$ is a cycle if
           and only if $\D'$ is a cycle.  
        \end{proof}

\begin{theorem}[Structure of a cycle]\label{t:cycle-structure}
 Let $\D$ be a facet complex.  Then $\D$ is a cycle if and only if
$\D$ can be written as a sequence of strong neighbors $F_1\sim_{\D}
F_2 \sim_{\D} \ldots \sim_{\D} F_n \sim_{\D} F_1$ such that $n\geq 3$,
and for all $i,j$ $$F_i \cap F_j = \bigcap_{k=1}^n F_k \quad\mbox{ if } j
\neq i-1, i,i+1\ (\mbox{\rm mod } n).$$
\end{theorem}

       \begin{proof} Let $\D$ be a cycle. Then by Proposition~\ref{p:cycle}
        and Lemma~\ref{lem:cycle-fund}, $\D$ can be written as a
        sequence of strong neighbors with the desired properties.

        Conversely, suppose that $\D$ is written as a sequence of
        strong neighbors $F_1\sim_{\D} F_2 \sim_{\D} \ldots
        \sim_{\D} F_n \sim_{\D} F_1$ such that $F_i \cap F_j =
        \bigcap_{k=1}^n F_k$ if $j \neq i-1,i, i+1$ (mod $n$).  By
        Lemma~\ref{l:cone} we can without loss of generality assume
        that $ \bigcap_{k=1}^n F_k=\emptyset$.

        By Corollary~\ref{c:cycle}, $\D$ has no leaf.  Suppose $\D'$ is
        any nonempty proper subset of $\D$. We need to show that $\D'$
        has a leaf. Suppose $F_i\in \D'$ and $F_{i+1}\not \in \D'$.
        There are two cases:

    \begin{enumerate}

       \item $F_{i-1} \notin \D'$. In this case, since $F_i \cap F_k
       =\emptyset$ for all $F_k \in \D'\srmv{F_i}$, $F_i$ is a leaf.

       \item $F_{i-1} \in \D'$. In this case, $F_i \cap F_k \subseteq
       F_{i-1}$ for all $F_k \in \D' \srmv{F_i}$, and so $F_i$ is
       again a leaf, this time with $F_{i-1}$ as a joint.
    \end{enumerate}

     So $\D$ is a cycle and we are done.
        \end{proof}

The implication of Theorem~\ref{t:cycle-structure} is that a
 simplicial cycle has a very intuitive structure: it is either a
 sequence of facets joined together 
 to form a circle in such a way that all intersections are pairwise
 disjoint (this is the case where the intersection of all the facets is the
 empty set in Theorem~\ref{t:cycle-structure}), or it is a cone over such
 a structure (Lemma~\ref{l:cone}).

\begin{example}\label{e:cycle-structure} The facet complex $\D$ is a cycle.
  The facet complex $\Gamma$ is a cycle and is also a cone over the
  cycle $\Gamma'$.

\begin{center}
$\D=\begin{tabular}{c}\epsfig{file=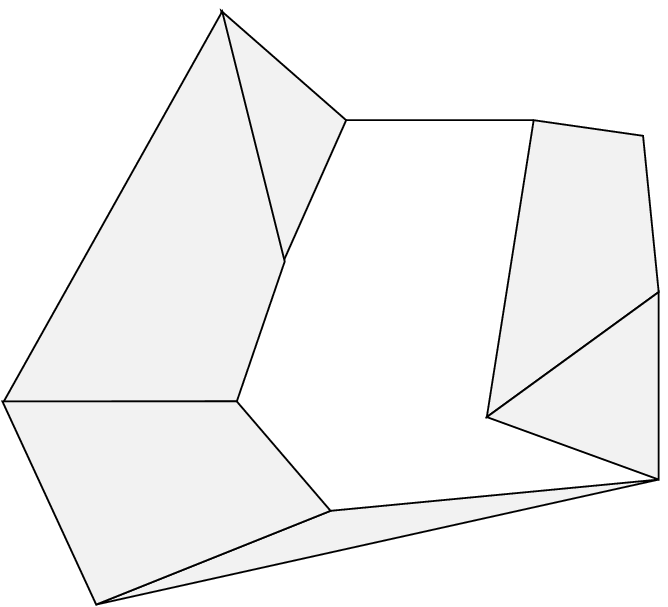, width=1.2in}\end{tabular}
\qquad
\Gamma=\begin{tabular}{c}\epsfig{file=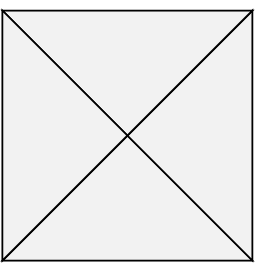, height=.7in}\end{tabular}
\qquad
\Gamma'=\begin{tabular}{c}\epsfig{file=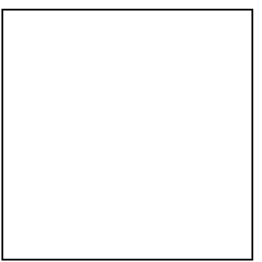, height=.7in}\end{tabular}$
\end{center}

\end{example}

\section{Characterization of trees}\label{s:tree-problem}

We now consider the problem of deciding whether or not a given facet
complex is a tree. We refer to this problem as the \emph{decision
  problem for simplicial trees}. 

Note that the na{\"\i}ve algorithm (namely, checking whether every
non-empty subset has a leaf) is extremely inefficient: for a facet complex
of $n$ facets, there are $2^n-1$ subsets to check. Also note that the
definition of a tree is not inductive in any obvious way: for
instance, attaching a single leaf to a tree need not yield a tree, as
Example~\ref{e:free-example} shows. This seems to rule out an easy
recursive algorithm.

Nevertheless, we demonstrate that the decision problem for simplicial
trees can be solved efficiently. This is done via a characterization
of trees given in this section. 

\begin{definition}[Paths and connectedness outside $V$]
  Let $\D$ be a facet complex, and let $V$ be a set of vertices.  We
  say that a sequence of facets $H_1,\ldots,H_n\in\D$ is a \emph{path
    outside $V$} in $\D$ if for all $i=1,\ldots,n-1$, $(H_i\cap
  H_{i+1})\sm V\neq\es$. We say that two facets $F,G\in\D$ are
  \emph{connected outside $V$} in $\D$ if there exists a path
  $H_1,\ldots,H_n$ outside $V$ in $\D$ such that $H_1=F$ and $H_n=G$.
\end{definition}

Note that in case $V=\es$, this coincides with the definition of
connectedness from Definition~\ref{def:connected}.

\begin{notation}
  If $F,G_1,G_2$ are three distinct facets of $\D$, then we define 
  $\residue{\D}{F}{G_1,G_2}$ to be the following subset of $\D$:
  \[ \residue{\D}{F}{G_1,G_2} = \s{H\in\D\such H\cap F=G_1\cap G_2}
\cup\s{G_1,G_2}. \]
\end{notation}

\begin{definition}[Triple condition]
  Let $\D$ be a facet complex. We say a triple of facets
  $\p{F,G_1,G_2}$ satisfies the \emph{triple condition} if
  $G_1\not\leq_F G_2$ and $G_2\not\leq_F G_1$, and if $G_1$ and $G_2$
  are connected outside $F$ in the facet complex
  $\residue{\D}{F}{G_1,G_2}$.
\end{definition}

We note that the definitions of $\residue{\D}{F}{G_1,G_2}$ and the
triple condition have changed from an earlier version of this article
{\cite{CFS}}; they have been simplified.

\begin{example} Consider $\D$ in Example~\ref{e:cycle-is}. Then the triple 
$\tuple{F_1,F_2,F_4}$ satisfies the triple condition. This is because
$F_4 \not \leq_{F_1} F_2$ and $F_2 \not \leq_{F_1} F_4$. Moreover
$\residue{\D}{F_1}{F_2,F_4}=\s{F_2,F_3,F_4,G}$, and a path connecting
$F_2$ and $F_4$ outside $F_1$ is $F_2,F_3,F_4$.

However, $\tuple{G,F_2,F_3}$ does not satisfy the triple condition,
since $F_2 \leq_G F_3$ (and $F_3 \leq_G F_2$). Also
$\residue{\D}{G}{F_2,F_3}=\s{F_2,F_3}$, and $F_2$ and $F_3$ are not
connected outside $G$.

\end{example}

\begin{proposition}[A triple is  part of a cycle]\label{p:triple-cycle}
  Let $\D$ be a facet complex. A triple $\p{F,G_1,G_2}$ satisfies the
  triple condition if and only if there exists a cycle $\D'\seq\D$
  such that $F,G_1,G_2\in\D'$ and $G_1\sim_{\D'}F\sim_{\D'}G_2$.
\end{proposition}

         \begin{proof} Suppose  $\p{F,G_1,G_2}$ satisfies the
          triple condition. Then by definition, $G_1 \not \leq_F G_2$
          and $G_2 \not \leq_F G_1$. Choose a minimal (with respect to
          inclusion) path
          $H_1,\ldots,H_n$ outside $F$ that connects $H_1=G_1$ to
          $H_n=G_2$. Note that minimality implies that for $j>i+1$, $(H_i
          \cap H_j) \sm F=\emptyset$. We claim
          that $\D'=\{F,H_1,\ldots,H_n\}$ is a cycle with
          \begin{equation}\label{e:chain} F\sim_{\D'}H_1 
          \sim_{\D'}\ldots \sim_{\D'}H_n
          \sim_{\D'}F.\end{equation}

          \begin{enumerate}
           \item[(a)] $F\sim_{\D'}G_1$ and $F\sim_{\D'}G_2$. 

                 If $F\cap G_1 \subseteq H_i$ for some $i$, $1<i<n$,
                 then since $H_i \in \residue{\D}{F}{G_1,G_2}$, we
                 have $F\cap G_1\seq H_i \cap F = G_1 \cap G_2\seq
                 G_2$. This implies that $G_1 \leq_F G_2$, a
                 contradiction. So $F\sim_{\D'}G_1$, and similarly
                 $F\sim_{\D'}G_2$

           \item[(b)] $H_i\sim_{\D'}H_{i+1}$ for $i=1,\ldots,n-1$. 
                 
                Since $(H_i \cap H_{i+1})\sm F\neq \emptyset$, we have
                that $H_i\cap H_{i+1} \not \seq F$.  By minimality of
                the path, if $H_i \cap H_{i+1} \subseteq H_j$ for some
                $j>i+1$, then $H_i\cap H_{i+1} \seq H_i\cap H_j \seq F$,
                a contradiction. The case $j<i$ is similar.
          \end{enumerate}
          This shows (\ref{e:chain}). To finish the proof that $\D'$
          is a cycle, we must show that it meets the remaining
          condition of Theorem~\ref{t:cycle-structure}. If $n=2$,
          there is nothing to show; assume therefore that $n\geq 3$.
           By definition of $\residue{\D}{F}{G_1,G_2}$, $F\cap
          H_j=G_1\cap G_2$ for $j=2,\ldots,n-1$, and so $\bigcap_{G\in
          \D'}G=G_1 \cap G_2$. Also, if $j>i+1$, then $H_i\cap H_j\seq
          F$ by minimality of the path, therefore
          $$H_i \cap H_j = (H_i\cap F) \cap (H_j\cap F)=G_1 \cap
          G_2=\textstyle\bigcap_{G\in
          \D'}G.$$
       So $\D'$ is a cycle.

      Conversely, suppose that $\D'$ is a cycle containing $F$, $G_1$
      and $G_2$, written as $F\sim_{\D'}G_1\sim_{\D'}H_1\sim_{\D'}
      \ldots \sim_{\D'}H_n\sim_{\D'}G_2 \sim_{\D'}F$, where $n\geq 0$.

       From the strong neighbor relations it follows that $G_1 \not
      \leq_F G_2$ and $G_2 \not \leq_F G_1$.  It also follows that the
      above sequence of strong neighbors provides a path from $G_1$ to
      $G_2$ outside $F$.  We only need to show that for
      $i=1,\ldots,n$, $H_i\cap F=G_1\cap G_2$.

      If $\D'=\{F,G_1,G_2\}$ we are done. So assume that $n\geq 1$.

      We know $H_i\not \sim_{\D'}F$, and so by
      Lemma~\ref{lem:cycle-fund}, $H_i\cap F\subseteq G_1 \cap
      G_2$. On the other hand, since $G_1 \not \sim_{\D'} G_2$,
      Lemma~\ref{lem:cycle-fund} implies the opposite inclusion
      $H_i\cap F\supseteq G_1 \cap G_2$. It therefore follows that
      $H_i\cap F = G_1 \cap G_2$ and we are done.
      \end{proof}

An immediate implication of Proposition~\ref{p:triple-cycle} is an
 (algorithmically) efficient criterion to determine whether or not a
 facet complex is a tree.

\begin{theorem}[Main Theorem]\label{t:tree-char}
  Let $\D$ be a connected facet complex. Then $\D$ is a tree if and
  only if no triple of facets in $\D$ satisfies the triple condition.
\end{theorem}

\section{A polynomial-time tree decision algorithm}\label{s:tree-algorithm}

By Theorem~\ref{t:tree-char}, to check if a facet complex
$\D=\{G_1,\ldots,G_l\}$ is a tree, we only need to check the triple
condition for all triples of elements of $\D$.  The checks themselves
are straightforward. Since the triple condition for $\p{F,G,G'}$ is
clearly unchanged if one switches $G$ and $G'$, we can limit triple
checking to the elements of the set $\{\p{F,G_i,G_j}\in \D^3 \mid
G_i\neq F\neq G_j,i<j\}$. The procedures for the basic steps follow
immediately from the earlier definitions.

\begin{algorithm}[Tree decision algorithm]\label{a:tree-algorithm}~\\
  Input: a connected facet complex $\D=\{G_1,\ldots,G_l\}$ with  $n$ vertices.\\
  Output: {\bf True} if $\D$ is a tree, {\bf False} otherwise.
  \begin{enumerate}
  \item For each triple $\p{F,G,G'}\in\{\p{F,G_i,G_j}\in \D^3 \mid
    G_i\neq F\neq G_j,i<j\}$
    \begin{enumerate}
    \item If $G\leq_F G'$ or $G'\leq_F G$, continue with the
      next triple.
    \item Build $\residue{\D}{F}{G,G'}$.
    \item If $G$ and $G'$ are connected outside $F$ in
      $\residue{\D}{F}{G,G'}$, return {\bf False}.
    \end{enumerate}
  \item Return {\bf True}.
  \end{enumerate}
\end{algorithm}

The correctness of this algorithm is an immediate consequence of
Theorem~\ref{t:tree-char}.  The algorithm uses very little memory; the
input $\D$ requires $nl$
bits, and $\residue{\D}{F}{G,G'}\seq\D$ requires $l$ bits.  The memory
required to perform the connectedness check 
and to store the various counters is negligible. Thus, memory
locality is good, and the computations can generally take place
in the cache.

\begin{remark}
  In the process of checking the triple condition for a triple
  $\p{F,G,G'}$ that is part of a cycle, we build a connection
  path outside $F$. Clearly, any such path can be reduced to a
  \emph{minimal} connection path $\s{H_1,\ldots,H_n}$ outside $F$ for
  $G,G'$, and therefore, by the proof of
  Proposition~\ref{p:triple-cycle}, $\{F,H_1,\ldots,H_n\}$ forms a
  cycle. Therefore, an easy modification of
  Algorithm~\ref{a:tree-algorithm} allow us to produce the set of all
  the facets $F\in\D$ that are part of some cycle, and a cycle
  $\D'_F\supseteq\{F\}$ for each of them.
\end{remark}


\subsection{Complexity}\label{s:tree-complexity} 

For each triple it is trivial to see that
steps (a) and (b) can be performed with cost $\OO{n}$ and  $\OO{nl}$
respectively. For step (c), the following holds.

\begin{lemma}\label{l:conn-check}
Let $\D$ be a facet complex with $l$ facets over $n$ variables such that
$F,G,G'$ are distinct facets of $\D$.
The connectedness outside $F$ of $G,G'\in\D$ can 
be determined with time cost $\OO{nl}$.
\end{lemma}
\begin{proof}
  First of all we substitute $\D$ with the set $\{H\sm F\mid H\in \D\}$.
  We then define $n+1$ equivalence relations $P_0,\ldots,P_n$ on
  the set $\s{1,\ldots,l}$. $P_0$ is the identity relation, i.e., each
  equivalence class is a singleton. For each $j=1,\ldots,n$, consider
  the vertex $v_j$ and the set $X_j=\s{i\such v_j\in F_i}$.  Let $P_j$
  be the smallest equivalence relation such that $P_{j-1}\seq P_j$ and
  such that for all $i,i'\in X_j$, $(i,i')\in P_j$. Then facets $F_i$
  and $F_{i'}$ are connected if and only if $(i,i')\in P_n$. With a suitable data
  structure for representing equivalence
  relations, the complexity of the procedure above is $\OO{nl}$.
\end{proof}

Consequently, step (c) of the tree decision algorithm can be performed
at cost $\OO{nl}$. Thus, the total complexity of the tree decision
algorithm is as follows: in the worst case we have to check $ 3\cdot
{l \choose 3}=\frac{l(l-1)(l-2)}{2}=\OO{l^3}$ triples. The complexity
of the steps (a)--(c) is $\OO{nl}$ and hence the total
complexity of the algorithm is $\OO{nl^4}$.

\begin{example}\label{e:triple}
Consider the facet complex $\D=\{xy,xz,yz,yu,zt\}$. We have to check
$3\cdot {5 \choose 3}=30$ triples. We start with the triple $\p{xy, xz, yz}$.
\begin{itemize}
\item $xz\not\leq_{xy}yz$ since 
$xy\cap xz=x\not\seq y =xy\cap yz$.
Similarly $yz\not\leq_{xy}xz$.
\item $xz$ and $yz$ are connected outside $xy$ in 
the complex $\residue{\D}{xy}{xz,yz}~=~\{zt,xz,yz\}$.
\end{itemize}
We have hence discovered that $\D$ is not a tree. A more unlucky
choice of facets could have brought about the checking of $27$ useless
triples before the discovery that $\D$ is not a tree, the other two
useful triples being $\p{yz,xy,xz}$ and $\p{xz,xy,yz}$.
\end{example}

\begin{example}\label{e:bigger-triple}   Some statistics for a bigger random example. Consider the
facet complex $\D=\{lka$, $qik$, $tykj$, $wuv$, $rjb$, $eioab$, $gdc$, $zv$,
$rtj$, $qrvm$, $gzm$, $tgzb$, $rgvm$, $qlav$, $qeocn$, $ikfaz$, $bn$,
$ekjs$, $pfvn$, $wtodv\}$. We discover that it is not a tree after
checking $4$ facets; we performed the connectedness check only
once. If one checks all $3\cdot {20\choose 3}=3420$ triples, one finds that
$445$ of them require a connectedness check, and $403$ of them reveal
that $\D$ is not a tree.
\end{example}

\begin{example}\label{line400}
The facet complex $\{x_ix_{i+1}x_{i+2}\mid i=1,\ldots, 400\}$ is trivially a
tree.  Checking this by a direct application of
Algorithm~\ref{a:tree-algorithm} requires dealing with
$3\cdot{400\choose 3}=31,760,400$ triples, and takes about $12.6$
seconds on an Athlon 2600+ machine for our C++ implementation. All the
timings in the remainder of this paper refer to this machine.
\end{example}


\subsection{Optimization}\label{sub:optimization} 

The runtime of Algorithm~\ref{a:tree-algorithm} can be improved by
introducing some optimizations. First, note that if $F$ is a facet
such that no triple $\p{F,G,G'}$ satisfies the triple condition, then
by Proposition~\ref{p:triple-cycle}, $F$ cannot be part of any
cycle of $\D$. Therefore, $F$ can be removed from $\D$,
reducing the number of subsequent triple checks. 
We refer to this optimization as the \emph{removal of useless facets}.

\begin{example}
We check the tree $\{x_ix_{i+1}x_{i+2}\mid i=1,\ldots, 400\}$ of Example~\ref{line400}
with a version of Algorithm~\ref{a:tree-algorithm} with removal of useless facets. This 
requires checking $10,586,800$ triples 
and takes about $3.46$ seconds. 
\end{example}

An important special case of a ``useless facet'' is a reducible leaf, as
captured in the following definition:

\begin{definition}[Reducible leaf]\label{d:reducible-leaf}
  A facet $F$ of a facet complex $\D$ is called a \emph{reducible
    leaf} if for all $G, G'\in\D$, either $G\leq_F G'$ or $G'\leq_F G$.
\end{definition}

  A reducible leaf is called a ``good leaf'' by Zheng~\cite{Z}.

\begin{remark}\label{r:reducible}
 The facet $F$ is a reducible leaf of $\D$ if and only if $F$ is a
  leaf of every $\D'\seq\D$ with $F\in \D'$.
\end{remark}

The remark immediately implies that a reducible leaf cannot be part of
a cycle. Thus, it can be removed from $\D$, and the algorithm can
then be recursively applied to $\D\rmf{F}$. 
We were not able to find a tree without a reducible leaf;
in fact, Zheng~\cite{Z} conjectured that this is always the case.  Checking
whether a given facet $F$ is a reducible leaf requires ordering all
facets with respect to $\leq_F$, which takes $\OO{nl\log l}$ steps. A
reducible leaf can thus be found in time $\OO{nl^2\log l}$. Therefore,
if Zheng's conjecture is true, the tree problem can be decided in time
$\OO{nl^3\log l}$. But even if the conjecture is not true, removing
all reducible leaves at the beginning of
Algorithm~\ref{a:tree-algorithm} is still a worthwhile optimization.

\subsection{Optimization for sparse complexes}\label{s:sparse-optimization} 

Let $\D$ be a facet complex with $l$ facets. If every $F\in \D$ intersects a
substantial ($\approx l$) number of facets, then the number of cycles is probably
high and our algorithm is usually able to detect one of them easily. If this
does not happen, we can exploit the ``sparseness'' of the facet complex in our
algorithm.

For the remainder of this subsection, $\D$ will be a facet complex
with $l$ facets over $n$ vertices such that 
the maximum number of neighbors of a facet $F\in\D$ is $d$ 
and the maximum number of vertices of a facet $F\in\D$ is $v$. 
Note that trees are the hard cases for our algorithm, since all the triples have
to be checked. 
Also note that, if $\D$ is a tree, then $l\leq n$. This follows by induction 
on $l$, from the fact that every leaf contains at least one free vertex.

\subsubsection{Connection set algorithm}
To check
if $\D$ is a tree it is sufficient to check the connected triples
only. For each facet $F$ ($l$ facets): first construct the set of all
facets $G$ connected to $F$ (called the \emph{connection set}, at cost
$\OO{lv}$), then for all $G,G'$ in the set ($\OO{d^2}$ pairs) perform the
triple check on $\p{F,G,G'}$ (cost $\OO{nl}$ per triple). We call this
optimization of Algorithm~\ref{a:tree-algorithm} the \emph{connection set algorithm}.
The total cost is $\OO{nl^2d^2}$. The space required to construct the connection
sets is $\OO{d}$, hence negligible.  If the complex is not sparse
($d\approx l$, $v\approx n$), the complexity is the same as 
Algorithm~\ref{a:tree-algorithm}. However, for sparse examples, this
optimization is clearly worthwhile:

\begin{example}
We check the tree $\{x_ix_{i+1}x_{i+2}\mid i=1,\ldots, 400\}$ of
Example~\ref{line400} with the algorithm detailed above.
We deal with $398$ triples and spend $0.2$ seconds.
\end{example}

\begin{example}   The facet complex  $\{x_ix_{i+1}\cdots x_{i+200}\mid
i=1,\ldots, 3200\}$ is a tree but not sparse. Tree checking with the connection
set algorithm is still quite efficient; it requires dealing with $61,013,400$
triples, and takes about $140$ seconds. Without any optimization, the number of
triples to check is $16,368,643,200$ and the time spent by the algorithm is
$>2$ days. 
\end{example}

\subsubsection{Incidence matrix algorithm}

The connectedness relation for a facet complex $\D$ can be represented by a
graph through an incidence matrix. This matrix can be built and used
during the tree checking algorithm. Since creating incidence matrices from
a complex is a relatively expensive operation, we build them in steps,
exploiting at every step the relations already computed.

We compute the connectedness relation for $\D$ at cost $\OO{l^2d}$. 
Then for every facet $F\in\D$
we compute the ``connectedness outside $F$'' relation for $\D$, at cost $\OO{nld}$.
Then for every triple $\p{F,G,G'}$ (there are $\OO{d^2}$ of them) we compute
the ``connectedness outside $F$'' relation for $\residue{\D}{F}{G,G'}$ at cost
$\OO{dv+ld}$.
Using this additional structure, we do not actually need to build 
$\residue{\D}{F}{G,G'}$, and we can check connectedness outside $F$ in 
$\residue{\D}{F}{G,G'}$ using the connectedness relations at cost $\OO{ld}$.
We call this optimization of Algorithm~\ref{a:tree-algorithm} the
\emph{incidence matrix algorithm}.

The total complexity for this algorithm is hence
$\OO{nl^2d+ld^3v+l^2d^3}$.  If $\D$ is not sparse ($v\approx n$,
$d\approx l$), then this algorithm has roughly the same complexity as
Algorithm~\ref{a:tree-algorithm}.

On the other hand, if $d\approx v\approx\sqrt{l}\approx\sqrt{n}$,
which is a reasonable assumption
for sparseness, then the complexity of the incidence matrix algorithm is
$\OO{l^3\sqrt{l}}$,  while the complexity of the connection set algorithm is
$\OO{l^4}$ and that of Algorithm~\ref{a:tree-algorithm}  is  $\OO{l^5}$.

\section{Algebraic properties of facet ideals}\label{s:CM-properties} 

We now study facet ideals from a more algebraic point of view. In
particular, we are interested in ways to determine whether a given
facet complex $\D$ is Cohen-Macaulay, meaning whether $R/\F(\D)$ is a
Cohen-Macaulay ring. We first need to introduce some new terminology.

\begin{definition}[Vertex covering number, unmixed facet complex] 
Let $\D$ be a facet complex. A \emph{vertex cover} for $\D$ is a set
  $A$ of vertices of $\D$, such that $A \cap F \neq \es$ for every
  facet $F$. The smallest cardinality of a vertex cover of $\D$ is
  called the \emph{vertex covering number} of $\D$ and is denoted
  by $\al(\D)$.  A vertex cover $A$ is \emph{minimal} if no proper
  subset of $A$ is a vertex cover. A facet complex $\D$ is
  \emph{unmixed} if all of its minimal vertex covers have the same
  cardinality.

\end{definition}

\begin{example} Consider the two facet complexes in Example~\ref{example11}.
We have $\al(\D)=2$. Also, $\D$ is unmixed as its minimal
vertex covers $\{x,u \}$, $\{y,u\}$, $\{y,v\}$, $\{z,u\}$ and
$\{z,v\}$ all have cardinality equal to two. We further have
$\al(\D')=1$, but $\D'$ is not unmixed, because $\{x\}$ and
$\{y,z\}$ are minimal vertex covers of different cardinalities.
\end{example}

The following observations are basic but useful.

\begin{proposition}[Cohen-Macaulay facet complexes~\cite{F1, F3}]\label{minprime}
Let $\D$ be a facet complex with vertices in $\xs$, and consider
its facet ideal $I=\F(\D)$ in the polynomial ring $R=k[\xs]$. Then
the following hold:
\begin{enumerate}

\item[(a)] $\height I\ = \al(\D)$ and $\dimn R/I =n -\al(\D)$.

\item[(b)] An ideal $p=(x_{i_1},\ldots,x_{i_s})$ of $R$ is a minimal
  prime of $I$ if and only if the set $\{x_{i_1},\ldots,x_{i_s}\}$ is
  a minimal vertex cover for $\D$.

\item[(c)] If $k[\xs]/\F(\D)$ is Cohen-Macaulay, then $\D$ is
unmixed.
\end{enumerate}
\end{proposition}

\subsection{Grafting}\label{s:grafting}

One of the most basic ways to build a Cohen-Macaulay facet complex is via
grafting.

\begin{definition}[Grafting \cite{F3}]\label{grafting} A facet  complex $\D$
  is a \emph{grafting} of the facet complex $\D'=\s{G_1,
    \ldots, G_s}$ with the facets $F_1, \ldots, F_r$ (or we say
    that $\D$ is \emph{grafted}) if $$\D= \s{F_1, \ldots, F_r}
    \cup \s{G_1, \ldots, G_s}$$ with the following properties:

\begin{enumerate}
\item[(i)] $G_1\cup\ldots\cup G_s \subseteq F_1 \cup  \ldots \cup F_r$;
\item[(ii)] $F_1, \ldots, F_r$ are all the leaves of $\D$;
\item[(iii)] $\{G_1, \ldots, G_s\} \cap \{F_1, \ldots, F_r\} = \emptyset$;
\item[(iv)] For $i \neq j$, $F_i \cap F_j = \emptyset$;
\item[(v)] If $G_i$ is a joint of $\D$, then $\D \rmf{G_i}$ is
  also grafted.

\end{enumerate} 
\end{definition}

Note that the definition is recursive, since graftedness of $\D$ is
defined in terms of graftedness of $\D\rmf{G_i}$. 
Also note that a facet complex that consists of only one facet or
several pairwise disjoint facets is grafted, as it can be
considered as a grafting of the empty facet complex. It is easy
to check that conditions (i) to (v) above are satisfied in this case.
It is also clear that the union of two or more grafted facet
complexes is itself grafted.

\begin{example}\label{grafting-example}  
There may be more than one way to graft a given facet complex.  For example,
  some possible ways of grafting $\s{G_1,G_2}$ are shown in Figure~\ref{f:grafting}.

\begin{figure}
\setlength{\unitlength}{0.28mm}
\begin{center}
\begin{picture}(385,220)
\put(0,180){$\D:$}
\put(40,155){\epsfig{file=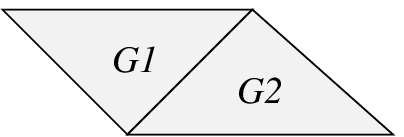, height=.4in}}
\put(150,180){$\stackrel{\mbox{\emph{graft}}}{\vector(1,0){50}}$}
\put(220,180){$\D':$}
\put(250,150){\epsfig{file=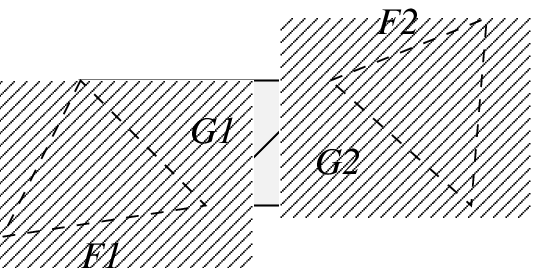, height=.8in}}
\put(100,130){$\vector(0,-1){50}$}
\put(110,105){\emph{graft}}
\put(0,30){$\D'':$} 
\put(40,0){\epsfig{file=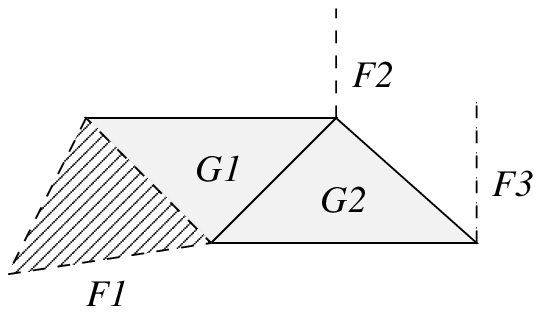, height=.85in}}
\put(180,130){$\vector(1,-1){50}$}
\put(220,105){\emph{graft}}
\put(220,30){$\D''':$}
\put(250,10){\epsfig{file=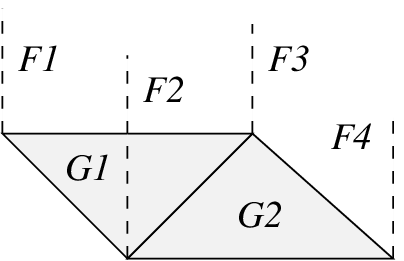, height=.8in}}
\end{picture}
\end{center}
\caption{Three different ways of grafting the facet complex $\D$.}\label{f:grafting}
\end{figure}

\end{example}

The interest in grafted facet complexes, from an algebraic point of view,
lies in the following facts.

\begin{theorem}[Grafted facet complexes are Cohen-Macaulay \cite{F3}]\label{th:CM} 
Let $\D$ be a grafted facet complex. Then $\F(\D)$ is
Cohen-Macaulay.\end{theorem}

Even more holds when $\D$ is a tree.

\begin{theorem}[{\cite[Corollaries~7.8,~8.3]{F3}}]\label{th:tree-unmixed-iff-grafted}
 If $\D$ is a simplicial tree, then the following are equivalent:
\begin{enumerate}
\item[(i)]$\D$ is unmixed;
\item[(ii)] $\D$ is grafted;
\item[(iii)] $\F(\D)$ is Cohen-Macaulay.
\end{enumerate}
\end{theorem}

\subsection{Graftedness  algorithm}\label{s:graft-algorithm}

A direct application of Definition~\ref{grafting} is not very
convenient for checking whether a given facet complex $\D$ is grafted,
since at each step of the recursion, one potentially needs to check
condition (v) for several of the $G_i$, and this leads to a worst-case
exponential algorithm. In order to arrive at a more efficient
algorithm, we characterize graftedness as follows:

\begin{lemma}[{cf. \cite[Remarks~7.2, 7.3]{F3}}]\label{l:grafted-equivalent}
  A facet complex $\D$ is grafted if and only if (1) for each vertex
  $v$, there exists a unique leaf $F$ such that $v\in F$, and (2) all
  leaves of $\D$ are reducible.
\end{lemma}

\begin{proof}[Sketch of the proof]
  First, assume that $\D$ is grafted. Condition (1) follows from (i),
  (ii) and (iv). The fact that all leaves are reducible is shown by
  induction on the number of facets of $\D$. The converse is also
  shown by induction. Suppose $\D$ satisfies (1) and (2), and let
  $\s{F_1,\ldots,F_r}$ and $\s{G_1,\ldots,G_s}$ be the sets of leaves
  and non-leaves, respectively. Conditions (i)--(iv) hold
  trivially. Further, if $G_i$ is a joint, then $F_1,\ldots,F_r$ are
  still reducible leaves of $\D\rmf{G_i}$ by Remark~\ref{r:reducible}.
  Also, there are no additional leaves in $\D\rmf{G_i}$, since none of
  the $G_j$ have free vertices by Condition (1). Therefore,
  $\D\rmf{G_i}$ satisfies (1) and (2) and is therefore grafted by
  induction hypothesis, proving (v).
\end{proof}

The algorithm for checking if a facet complex is grafted follows immediately
from Lemma~\ref{l:grafted-equivalent}.

\begin{algorithm}[Graftedness algorithm]\label{a:grafted-algorithm}
  ~\\ Input: A facet complex $\D$ with $l$ facets and $n$ vertices.\\
  Output: {\bf True} if $\D$ is grafted, {\bf False} otherwise.
  \begin{enumerate}
  \item Build the lists $\calF=\s{F_1,\ldots,F_k}$ (leaves of $\D$) and 
    $\calG=\s{G_1,\ldots,G_m}$ (facets of $\D$ which are not leaves).
  \item If $\bigcup_{G\in\calG} G\not\seq\bigcup_{F\in\calF} F$, return {\bf False}.
  \item If $\exists\ F,F'\in\calF\text{ such that } F\cap F'\neq\emptyset$, return {\bf False}.
  \item If $\exists\ F\in\calF$ that is not a reducible leaf, return
  {\bf False}.
  \item Return {\bf True}.
  \end{enumerate}
\end{algorithm}

\subsection{Complexity}\label{s:graft-complexity} 

The leaf checking cost is $\OO{nl}$, hence the cost of step 1 is
$\OO{nl^2}$. The cost of steps 2 and 3 is $\OO{nl}$. For step 4, there
are $k$ facets $F$ to check. Checking whether $F$ is reducible takes
$\OO{nl\log l}$ steps as mentioned in
Section~\ref{sub:optimization}. Therefore the total cost for step 4 is
$\OO{nl^2\log l}$, and this is the cost of the algorithm.

\begin{example} Let $\D=\{xyz,yzu,ztu,uv,tw\}$,
with $\calF=\{xyz, uv, tw\}$ and $\calG=\{yzu, ztu\}$.  Then
 $\bigcup_{G\in\calG} G \subseteq \bigcup_{F\in\calF} F= \{x, y, z, t,
 u, v, w\}$ and $xyz\cap uv=xyz\cap tw=uv\cap tw=\emptyset$.
 Additionally, we check that each $F \in \calF$ is a reducible leaf by
 showing that the set $\{F \cap G \st G \in \calG\}$ is a totally
 ordered set under inclusion. For example, if $F=xyz$, then this set
 is equal to $\s{yz,z}$ which is totally ordered. This holds for all
 $F \in \calF$, and hence the facet complex is grafted.
 \[
 \epsfig{file=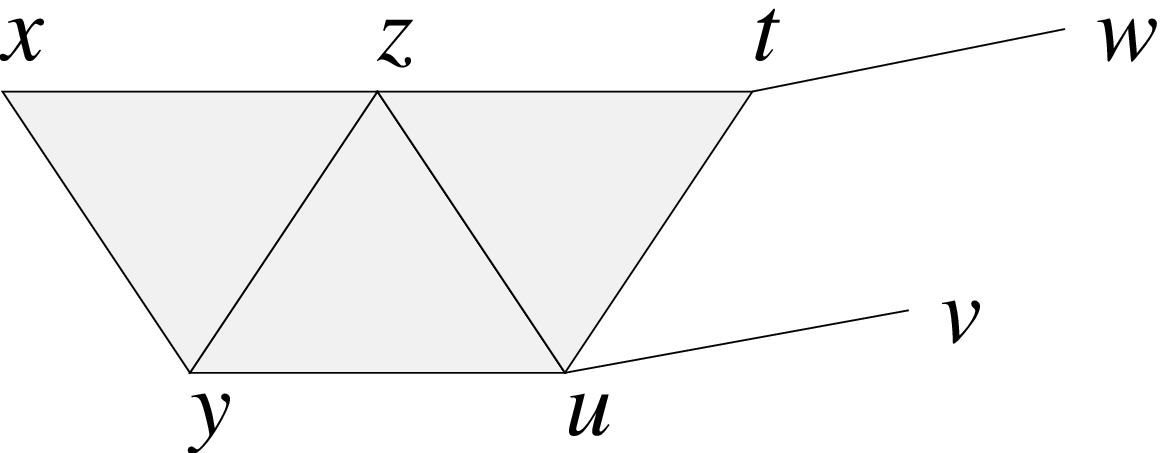, height=.6in}
 \]
\end{example}



\end{document}